\documentclass[11pt]{article}

\usepackage{amsmath}
\usepackage{amssymb}
\usepackage{amscd}

\newtheorem{teor}{Theorem}[section]

\newtheorem{lemma}[teor]{Lemma}

\newcommand{\ann}{\operatorname{Ann}}
\newcommand{\G}{\mathbb{G}}
\newcommand{\Gr}{{\G(r,N)}}
\newcommand{\vol}{\operatorname{vol}}

\newcommand{\Hom}{\operatorname{Hom}}

\renewcommand{\dim}{\operatorname{dim}}

\newcommand{\SU} {\operatorname{SU}}
\newcommand{\su} {\mathfrak{su}}

\newcommand{\SL}{\operatorname{SL}}

\newcommand{\fdim}{\hspace*{\fill}$\Box$}
\newcommand{\dimostr}[1][]    {\noindent\textbf{Proof#1}. }
 
\newcommand{\real}{\mathbb{R}}

\newcommand{\complex}{\mathbb{C}}
\newcommand{\projective}{\mathbb{P}}

\newcommand{\K}{K\"{a}hler}
\newcommand{\GL}{\operatorname{GL}}
\newcommand{\tr}{\operatorname{tr}}
\newcommand{\ra}{\rightarrow}
\newcommand{\C}{\mathbb{C}}
\newcommand{\debar}{\bar{\partial}}
\newcommand{\barz}{\bar{z}}
\newcommand{\OO}{\mathcal{O}}
\newcommand{\unders}{{\underline{s}}}
\newcommand{\KE}{K\"ahler-Einstein }
\newcommand{\om}{\omega}
\newcommand{\eps}{\varepsilon}
\renewcommand{\phi}{\varphi}
\newcommand{\W}{\mathbb{W}}
\newcommand{\Z}{\mathbb{Z}}
\newcommand{\scalar}{\langle \cdot \, , \cdot \rangle}

\newcommand{\fine}    {\begin{flushright}
             \textsc{Q.E.D.}
             \end{flushright}}

\newcommand{\basis}[1]{{{\boldsymbol{#1}}}}

\begin{document}

\title{Stable bundles and the first eigenvalue of the Laplacian}

\author{Claudio Arezzo, Alessandro Ghigi, Andrea Loi} \date{}
\maketitle

{
\abstract{
\noindent In this paper we study the first eigenvalue of the Laplacian on a 
compact  manifold using stable bundles
and balanced bases. Our main result is the following: let $M$ be a 
compact \K\ manifold of complex dimension $n$ and $E$ a holomorphic 
vector bundle of rank $r$ over $M$. If $E$ is globally generated and 
its Gieseker point $T_E$ is stable, then for any \K\ metric $g$ on $M$
\begin{equation*}
\lambda_1(M, g) \leq
\frac{ 4  \,  \pi \, \, h^0(E)  }     {r(h^0(E)-r)}
\cdot
\frac{ \langle c_1(E) \cup      [\omega]^{n-1} , [M]\rangle}
      {(n-1)!  \,\, \vol(M, [\omega])},
    \end{equation*}
where $\omega =\omega_g$ is the \K\ form associated to $g$.

By this method we obtain, for example, a sharp upper bound for
$\lambda_1$ of \K\ metrics on  complex Grassmannians.

  }}

\section{Introduction and statements of the results}

The first eigenvalue of the Laplace operator is one of the most
natural and studied riemannian invariants.
A general question is how the underlying differentiable and topological
structure is sensitive of this riemannian invariant.
For example, in the case of a compact surface $M^2$, it is known that the
product
$\lambda_1(M,g) \cdot \vol(M,g)$ is bounded above by a function of the
genus only \cite{yang-yau}.
While for general riemannian metrics on higher dimensional manifolds
similar results
(substituting the genus with any geometric quantity) cannot hold, as proved
in  \cite{colbois-dodziuk},
it is natural to try to get upper bounds for some restricted natural
classes of metrics.

When the underlying compact manifold is K\"{a}hler one studies the
functional $\lambda_1(M,g) \cdot \vol(M,g)$ within a fixed \K\ class
$\alpha$ allowing also complex invariants to appear in the upper
bound. Since the volume is constant, this amounts to study the
behaviour of $\lambda_1$ on the metrics $g$ with $\om_g \in \alpha$.
%
%
%
%
By the Rayleigh principle, upper bounds for the first eigenvalue are
obtained by constructing functions with zero mean, sensitive to the
geometry of the underlying manifold.  If the manifold admits a map to
a space where one is able to produce abundance of mean zero functions,
one can try to import them on the manifold, thus getting an upper
bound. This idea has been the core of the work of Hersch
\cite{he} and Yang-Yau \cite{yang-yau}
for Riemann surfaces, using holomorphic maps onto $S^2$.  This strategy
ultimately relies on the possibility of using conformal
diffeomorphisms of the two sphere without altering the holomorphicity
of the map, in order to cover {\em all} riemannian metrics on the
Riemann surface.

It is clear that this strategy needs a serious change for higher
dimensional base manifolds.  Such a generalization is possible if one
restricts to \K\ metrics and the manifold is immersed in a projective
space. Bourguignon, Li and Yau \cite{bourguignon-li-yau} have in fact
shown that one can move the algebraic manifold with an automorphism of
the projective space in such a way that the pull-back of a certain
family of functions on $\projective{}^N$ have zero mean.  These
functions are in fact the components of the moment map for the action
of $\SU(N+1)$ on $\projective{}^N$.  It follows that for metrics $g$
with $\om_g\in\alpha$ the first eigenvalue is bounded above by an
invariant depending on the immersion in projective space.

In the light of recent results of Xiaowei Wang
\cite{wang-xiaowei-balance}, we observed that the theorem of
Bourguignon, Li and Yau can be rephrased as a suitable stability
property (stability of the Gieseker point) of any ample globally
generated line bundle.  This notion becomes in fact more interesting
for higher rank vector bundles and the aim of this paper is precisely
to investigate how we can use Gieseker stable vector bundles on \K\
manifolds to improve the upper bounds on $\lambda_1$.

First observe that vector bundles give rise to maps (if globally
generated) to Grassmannians.
Roughly speaking Wang proved that the Gieseker point of a globally
generated vector bundle is stable if and only if the associated map
into the Grassmannian can be moved into a ``balanced" position. Once
this is achieved we can give the seeked upper bound:
\begin{teor}\label{mainteor}
  Let $E\rightarrow M$ be a holomorphic vector bundle of rank $r$ over
  a compact \K\ manifold $M$ of complex dimension $n$. Assume that
  \begin{enumerate}
  \item $E$ is  globally generated,
  \item the Gieseker point $T_E$ is stable.
  \end{enumerate}
  Then, for any \K\ metric $g$ on $M$ one has
  the eigenvalue estimate
  \begin{equation}
    \label{mainest}
\begin{gathered}
    \lambda_1(M, g) \leq
    \frac{ 4\pi  \, h^0(E)  }     {r(h^0(E)-r)}
    \cdot
    \frac{ \langle c_1(E) \cup      [\omega]^{n-1} , [M]\rangle}
    {(n-1)!\vol(M, [\omega])}\,\, ,
  \end{gathered}
\end{equation}
where $\omega =\omega_g$ is the \K\	 form associated to $g$.

In particular, if $\omega_g \in 2\pi c_1(L)$ for some line bundle $L$ over
$M$, then
\begin{equation}
    \label{mainest2}
\begin{gathered}
    \lambda_1(M, g) \leq
  \frac{2n\, h^0(E) \deg E}{r(h^0(E)-r)c_1(L)^n} \,\, ,
  \end{gathered}
\end{equation}
where $\deg E=c_1(E) \cdot c_1(L)^{n-1}$.
\end{teor}
If $E$ is a line bundle this result reduces to the following
generalized version of Bourguignon-Li-Yau's estimate:
\begin{teor}\label{bly}(Bourguignon--Li--Yau)
  Let $g$ be a \K\ metric on a \K\ manifold
  $M$ and let $E$ be a globally
  generated line bundle on $M$ with \linebreak $N=h^0(E)=\dim H^0(E)$.  Let
  $\phi_\basis{t}:M\rightarrow {\complex}P^{N-1}$ be the Kodaira map
  in a basis $\basis{t}=(t_1,\dots , t_N)$ of $H^0(E)$.  Then
  \begin{equation}
    \label{lambda1yau}
    \lambda_1 (M, g)\leq \frac{4nN}{N-1}d.
  \end{equation}
Here $d$ is the so-called {\em holomorphic immersion degree} defined
by
$$d=
\frac{\int_M \phi_{\basis{t}}^*(\sigma)\wedge \omega^{n-1}}{\int_M
  \omega^n},$$
$\omega =\omega_g$ and $\sigma$ is half of the Fubini--Study form

$$\omega_{FS}=i\partial\bar\partial\log (|z_0|^2+\cdots +|z_{n-1}|^2)$$
on ${\complex}P^{N-1}$
(and hence $[\sigma]=\pi c_1(\OO(1))$).
If moreover $\omega\in 2\pi c_1(L)$ as above, then $d = \frac{\deg E}{2\,
c_1(L)^n}.$
\end{teor}

If $M=\projective{}^n$ and $L$ is the hyperplane bundle the estimate
given by Theorem \ref{bly} is sharp, since the bound is realized by
the Fubini-Study metrics.

It is then natural to suspect that our result gives a sharp estimate
for the Grassmannian manifolds, which are now the natural target
manifolds. This is indeed the case:
\begin{teor}\label{teorgrass}
  For any \K\ form $\omega_g$ on $M=\G (r, N)$ in the class $2\pi c_1
  (M)$  one has:
  $$\lambda_1(M, g)\leq 2.$$
\end{teor}
It is easy to see that this bound cannot be achieved using line
bundles and the theorem of Bourguignon, Li and Yau.  Moreover the
value $\lambda_1=2$ is indeed achieved by the symmetric \KE metric on
$\G(r,N)$. Thus the symmetric metric is a maximum point of the
functional $\lambda_1$ restricted to the set of \K\ metrics with fixed
volume.  This should be compared with a result of El Soufi and Ilias
according to which the symmetric metric is a critical point (in
suitable sense) for $\lambda_1$ on the set of \emph{all} Riemannian
metrics with fixed volume  (see Remark 1 at p. 96 of
\cite{el-soufi-ilias-Pacific}).

As mentioned above, our results are obtained by using balanced maps
into Grassmannians, hence implicitely using en passant a special type
of metrics on $M$ (those obtained as pull back
of the symmetric metric via these embeddings).
The main point
regarding these metrics is their abundance, since, as we prove
adapting Wang's work, they
are sensitive only of the stability of the Gieseker point of the vector
bundle.
This allows to prove general
estimates for {\em all} \K\ metrics on $M$. On the other hand a stronger
notion of ``balanced"
metric has been deeply studied in the last few years for line bundles
(see e.g. \cite{do}, \cite{al}, \cite{luo} and references therein). While
in general more rare
(they in fact measure a non vacuous stability property of line bundles),
it would be very interesting to know whether these metrics, when they
exist, can give stronger informations on the
spectrum of the laplacian. We leave this, and other question discussed
in the last section, for future research.

\paragraph{Acknowledgemts:}
 We wish to thank Gian Pietro Pirola for many enlightening discussions
concerning various aspects of this work. The first author wishes to thank also 
Kieran O'Grady for many helpful discussions about stability constructions.

\section{The set up and the proofs}

\subsection{Grassmannians}

\noindent If $W$ is a complex vector space we denote by $\G(r,W)$ the
Grassmannian of $r$-dimensional \emph{subspaces} of $W$. When $W=\C^N$
we will write $\G(r,N)$ for $\G(r,\C^N)$.  We denote by $U=U_{r,N} \ra
\G(r,N)$ the \emph{universal subbundle}. It is the subbundle of the
trivial bundle $\Gr\times \C^N$ whose fibre over a point $x\in \Gr$ is
simply the subspace $U_x=x$ represented by $x$.  If $a_1, \ldots ,
a_r$ is a basis of $U_x$, consider the $N\times r$ matrix $A=(a_1,
\ldots , a_r)$, whose columns are the vectors $a_\alpha$. This means
that if $a_\alpha=(a_{1\alpha},\ldots, a_{N\alpha})$ are the
components of the vector $a_\alpha$, then $A=(a_{i\alpha})$. (We let
the Greek indices run over $1,...,r$ and the Latin ones over $1,\ldots
, N$.)  We say that $A(x)$ is a \emph{Stiefel matrix} for the point
$x$ or that $A(x)$ are \emph{Stiefel coordinates} for $x$ (see
\cite{gelfand-kapranov-zelevinsky}).  If $a'_1, ..., a'_r$ is another
basis of $U_x$, there is a nonsingular matrix $C=(c_{\alpha\beta})\in
\GL(r,\C)$, such that $a'_\beta= c_{\alpha\beta}a_\alpha$. This means
that $A'=AC$. Therefore the Stiefel coordinates are defined only up to
right multiplication by a nonsigular $r\times r$ matrix.  This simply
reflects the fact that $\Gr= M^*(N,r,\C)/ \GL(r,\C)$, where
$M^*(N,r,\C)$ denotes the set of $N\times r$ matrices of maximal rank.
In terms of Stiefel coordinates the standard action of $ \GL(N,\C)$ on
$\Gr$ reads as follows: let $x$ be a point in $\Gr$, $P$ an element of
$\GL(N,\C)$ and $A$ a Stiefel matrix for $x$; then $PA$ is a Stiefel
matrix for $Px$.

Let now $e_1,\ldots , e_N$ be the standard basis of $\C^N$. For $I$ a
multiindex of lenght $r$ put
\begin{equation*}
  U_I=\{x\in \G(r,N): U_x + \mathrm{span} ( e_i : i\notin I ) = \C^N \}.
\end{equation*}
For simplicity of notation assume $I=(0,...,0)$. If $x\in U_0$ let $A$
be some Stiefel coordinates of $x$. Then
$$
A= \begin{pmatrix}A_1 \\ A_2
\end{pmatrix}
$$
where $A_1\in \GL(r,\C)$ and $A_2$ is an $(N-r)\times r$ matrix. The
affine coordinates of $x$ in the chart $U_0$ are then given by the
matrix $Z=A_1^{-1}A_2$. In other words the affine coordinates are the
entries of the matrix $Z$ such that $(I_r \ Z^t)^t$ be a Stiefel
coordinate of $x$.  Of course affine coordinates are honest
holomorphic coordinates for the complex manifold $\Gr$.

If $x\in \Gr$ and $a_1,\ldots, a_r$ is a basis of $U_x$, then
$a_1\wedge \cdots \wedge a_r$ is a nonzero element of $\Lambda^r\C^N$
defined up to multiplication by a nonzero scalar. Therefore it
represents a well-defined point in $\projective{}\bigl ( \Lambda^r
\C^N\bigr )$.  The corresponing map is an embedding of $\Gr$ in this
projective space, and it is called the \emph{ Pl\"ucker embedding}.
Observe that if we let $\OO(1)$ be the hyperplane bundle on
$\projective{}\bigl ( \Lambda^r \C^N\bigr )$ then $\det U^* =
\OO_\Gr(1)$ and $K_\Gr = \OO_\Gr(-N)$.

The constant metric on the fibres of $\Gr\times \C^N$ induces a
Hermitian metric on the subbundle $U $. Let $H_\G$ and $h_\G=\det
H_\G$ be the induced metric on $U^*$ and $\det U^*=\OO_\Gr(1)$
respectively.  Put $\omega_\G = iR(h_\G)$.  Let $A=(a_{ i\alpha})$ be
a Stiefel matrix for $x\in \G(r,N)$ and let $a_\alpha$ be the columns
of $A$. Then
\begin{equation*}
  ||a_1\wedge \cdots \wedge a_r ||_{h^*_\G}^2 = \det (A^*A),
\end{equation*}
where $h_{\G}^*$ is the induced metric on $\det U$.  Assume for
simplicity of notation that $x\in U_0$, and let $Z$ be the affine
coordinates of $x$. Then we can choose $A=(I_r\ Z^t)^t$. The
corresponding basis of $U_x$ is $\{e_\alpha + z_\alpha\}$, where
$z_\alpha$ is the $\alpha$-th column of $Z$. Then $s=(e_1+z_1)\wedge
\cdots \wedge (e_r+ z_r)$ is a nonzero section of $\Lambda^rU$ over
$U_0$ and $||s||^2 = \det (I_r + Z ^*Z)$. Therefore
\begin{equation*}
  \omega_\G = i \partial\debar \log \det (I_r + Z^*Z).
\end{equation*}
At the point $x=(I_r \ 0)^t$ the expression of $\omega_\G$ in the
affine coordinates $Z=(z_{p\alpha})$ is simply
\begin{equation}
  \label{eq:FS-affine}
  \omega_\G(x) = i \sum_{p=1}^{N-r}\sum_{\alpha=1}^r d z_{p\alpha}\wedge
  d \barz _{p\alpha}.
\end{equation}
When $r=1$, $\omega_\G$ is just the Fubini-Study metric.  When $r>1$
it is the pull-back of the Fubini-Study metric on
$\projective{}^{\binom{N}{r}-1}$ via the Pl\"ucker embedding.

The standard action of $\SU(N)$ on $\Gr$ is holomorphic and preserves
$\omega_\G$. Its moment map $ \mu_\G : \Gr \ra \su(N) $ is given by
\begin{equation}
  \label{mu-G}
  \mu_\G (x) = i \Bigl ( A(A^*A)^{-1}A^* - \frac{r}{N} I_N\Bigr ).
\end{equation}
Here $A$ are Stiefel coordinates of $x$ and we identify $\su(N)^*$
with $\su(N)$ by means of the Killing scalar product $\langle X,
Y\rangle = \tr X^*Y = - \tr XY$.  The normalization in chosen so that
\begin{equation}
  \label{eq:medianulla}
  \int_\Gr \mu_\G\operatorname{\vol}_\G =0.
\end{equation}
We use the sign convention so that
\begin{equation}
  \label{eq:moment}
  d\langle \mu,v\rangle = - i_{\xi_v}\omega_\G
\end{equation}
where $v\in \su(N)$ and $\xi_v$ is the fundamental vector field of the
action on $\Gr$.

For later use we also recall also the following identity:
\begin{equation}
  \label{omG-mu}
  \omega_\G = - i \sum_{j,k=1}^N d\mu_{jk}\wedge d \mu_{kj}.
\end{equation}
where $\mu_{jk}$ are the entries of $\mu_\G$.  Just as for \eqref
{eq:FS-affine} and \eqref{mu-G}, it is enough to prove this formula at
one point, thanks to the equivariance properties.  A simple
computation in affine coordinates shows that it holds at the origin of
the chart.

\subsection{Kempf-Ness theorem}

Let $G$ be complex reductive group, $K\subset G$ a maximal compact
subgroup, $W$ a linear representation of $G$ and $\langle \, , \,
\rangle$ a $K$-invariant Hermitian product on $W$.  For $v\in W$ the
function $\rho_v(g)=\log || g^{-1} v||$ is $K$-invariant and descends
to a convex function $\nu_v$ on the symmetric space $X=G/K$.  The moment
map $\mu: \mathbb{P}(W)\ra \mathfrak{k}$ for the action of $K$ is
given by $\mu([v]) = (d \rho_v)_e =(d \nu_v)_e$ where we consider $(d
\nu_v)_e$ as an element of $(\sqrt{-1}\mathfrak{k})^* = \mathfrak{k}$.

A point $x=[v]\in \projective{}(W)$ is said to be \emph{semistable} if
$0 \notin\overline{G.v} $, and it is called \emph{stable} if the orbit
$G\cdot v$ is closed and the stabiliser $G_v$ is finite.
\begin{teor}\label{kempf-ness}
  (Kempf-Ness) A point $x=[v]\in \projective(W)$ is semistable if and
  only if the function $\nu_v$ is bounded below. It is stable if and
  only if $\nu_v$ is proper if and only if $\mu$ has a unique zero on
  $G\cdot x$.

\end{teor}

\subsection{Vector bundles, stability and balanced bases}

Let $E$ be a holomorphic vector bundle of rank $r$ over a compact
complex manifold $M$.
Let $V=H^0(E)$ be the space of global holomorphic sections of $E$.
Assuming that $E$ is globally generated, for each $x\in M$ the
subspace $V_x\subset V$ of sections vanishing at $x$ is an
$(N-r)$-dimensional subspace. Denote by $\ann(V_x)$ its annihilator,
that is
\begin{equation*}
  \ann(V_x)=\{\lambda \in V^* : \lambda \equiv 0 \text { on } V_x\}.
\end{equation*}
Then one can define the so-called {\em Kodaira map}
\begin{equation*}\label{iv}
  \phi_{V}:M\rightarrow \G (r, V^*),\quad
  x\mapsto \ann(V_x).
\end{equation*}
A choice of a basis $s_1,\dots s_N$ of $V$ identifies $\G (r, V^*)$
with $\G (r, N)=\G(r, {\complex}^N)$. We denote by
$$
\phi_{\basis{s}}:M\rightarrow \G (r, N)
$$
the map $\phi_{V}$ written in the basis $\basis{s}=(s_1,\dots s_N)$.
If $\basis{\sigma}=(\sigma_1,\dots ,\sigma_r)$ is a local frame
for $E$ on a trivializing open set $U_{\basis{\sigma}}\subset M$
one has:
\begin{equation}
\label{trivializing}
s_j=\sum_{\alpha =1}^ra_{j\alpha}\sigma_{\alpha},\ j=,1\dots, N.
\end{equation}
If $A(x)$ denotes the $N\times r$ matrix with complex entries
$a_{j\alpha }(x)$ then $A: U_{\basis{\sigma}} \ra M^*(N,r,\C)$ is
a local expression of $\phi_{\basis{s}}$ in the Stiefel coordinates,
that is $A(x)$ is a Stiefel matrix for the point
$\phi_{\basis{s}}(x)$.  By construction we have the commutative
diagram
\begin{equation*}
  \begin{CD}
    E=\phi_{\basis{s}}^*U_{r, N}^* @>>> U_{r,N}^*\\
    @VVV @VVV\\
    M @>\phi_{\basis{s}}>>\G(r,N).
  \end{CD}
\end{equation*}
Denote by $k_{\basis{s}}$ the pull-back of $H_\G$ via the map
$\phi_{\basis{s}}$ namely
\begin{equation}\label{kFS}
  k_{\basis{s}}=\varphi_{\basis{s}}^*(H_{\G}).
\end{equation}
Let $\omega =\omega_g$  be a fixed \K\ form on $M$.  A basis
$\basis{s}=(s_1,\dots s_N)$ of $V$ is called
$\omega$-\emph{balanced} iff
\begin{equation}\label{oomegabal}
  \langle s_j , s_k\rangle_{k_{\basis{s}}, \omega}:=
  \int_M k_{\basis{s}}
  (s_j, s_k)\frac{\omega^n}{n!}=
  \frac{r\vol (M, g )}{N}\delta_{jk}.
\end{equation}
In words a basis $\basis{s}$ of $V$ is $\omega$-{\em balanced}
iff, up to the product with the positive constant $(r/N)\vol(M,g)$, 
it is an orthonormal basis of $V$ with respect to the
$L^2$-product $\langle\cdot , \cdot \rangle_{k_{\basis{s}},
  \omega}$ defined using the metric $k_\unders$ and the volume
$\omega^n/n!$. 

If $\sigma_1, \ldots{} , \sigma_r$ is a local frame for $E$ and $s_j$
are given by \eqref{trivializing} then
\begin{equation*}
\Bigl(  k_{\basis{s}}\bigl( s_j(x), s_k(x) \bigr )
\Bigr )_{j,k} =
A (AA^*)^{-1} A^*
\end{equation*}
that is
$$
\mu_{jk}\bigl( \phi_{\basis{s}}(x) \bigr ) = k_{\basis{s}}\bigl(
s_j(x), s_k(x) \bigr ) - \frac{r}{N}\,\delta_{jk}.
$$
Therefore the basis $\basis {s}$ is $\om$-balanced if and only if
\begin{equation}
  \label{eq:balanced-moment}
  \int_M \mu_\G \circ \phi_{\basis{s}} \frac{\om^n}{n!} = 0.
\end{equation}

An important theorem of Xiaowei Wang (conjectured by Donaldson in the
case of Riemann surfaces) relates $\omega$-balanced metrics to stable
bundles when $\omega \in 2\pi c_1(L)$.  In order to state it (see Theorem \ref{wuang0} below)
we recall some
definitions.  Let
$(M,L)$ be a polarised projective manifold. This means that $M$ is a
K\"ahler manifold and $L$ is an ample line bundle over $M$.  Given a
holomorphic vector bundle $E$ over $M$ of rank $r$, define the
\emph{degree} $\deg E$ and the \emph{slope} $\mu(E)$ of $E$ (with
respect to the polarisation $L$) by the formulas
\begin{equation}
  \deg E = c_1(E) \cdot c_1(L)^{n-1} \qquad \mu(E) = \frac{\deg E}{r}.
\end{equation}
Set $E(m)= E\otimes L^m$ and $ p_E(m) = (1/r)\cdot \chi(E(m)) $, where
$\chi$ denotes the Euler characteristic of a sheaf.  The vector bundle
$E$ is said to be \emph{Gieseker stable} if for any coherent subsheaf
$F \subset E$, and $m$ sufficiently large (depending on $F$) one has
the inequality $p_F(m) < p_E(m) $. On the other hand $E$ is said to be
\emph{Mumford-Takemoto stable } (or simply \emph{Mumford stable} or
\emph{slope stable}) if for any coherent subsheaf $F$ of $E$ one has
$\mu(F) < \mu( E)$. Since $\mu( E)$ is the leading coefficient of
$p_E$, Mumford stability is a stronger condition than Gieseker
stability (see e.g. \cite{friedman-universitext}, Chapter 4).  In
order to study the Gieseker stability of a globally generated bundle
$E$, Gieseker \cite{gieseker-vector-surfaces} considered the linear map
\begin{equation}
  T_E: \Lambda^r H^0(E) \longrightarrow H^0( \det E), \ 
(s_1, \dots ,s_r) \mapsto s_1 \wedge \cdots \wedge s_r.
 \end{equation}
This map, regarded as a point in $\mathbb{P}(\Hom \bigl (\Lambda^r
H^0(E), H^0(\det E)\bigr )$ is called the \emph{Gieseker point} of
$E$. On this projective space there is a natural action of
$\SL\bigl(V\bigr), V=H^0(E)$, and so it makes
sense to speak about the stability of $T_E$. 
Observe that the  stability of the Gieseker point does not involve the choice of an
ample line bundle $L$.

In the case where $X$ is a projective surface, Gieseker showed that a
vector bundle on $(X,L)$ is Gieseker stable if and only if $T_{E(m)}$ 
is stable for sufficiently large $m$ (see Theorem 0.7 in
\cite{gieseker-vector-surfaces}).  Wang, on the other hand, used the
Gieseker point to prove the following theorem (see Theorem 1.1 in
\cite{wang-xiaowei-balance}).

\begin{teor}\label{wuang0}(Wang)
Let $E$ be a holomorphic vector bundle on a polarised
projective manifold $(M, L)$ and let $\omega\in 2\pi c_1(L)$ be a \K\ form. 
Then $E$ is Gieseker stable iff there is a $m_0$
such that for all $m\geq m_0$ $E(m)$
admits an $\omega$-balanced basis.
\end{teor}

The result we actually need is the following 
slightly different version of this theorem,
where the polarization does not play any role.

\begin{lemma}\label{wang-teo} Let $E$ be a holomorphic vector
 bundle of rank $r$ over a compact \K\ manifold $M$, and let
  $\omega$ be a \K\ form on $M$.  If
  $E$ is globally generated and 
  the Gieseker point $T_E$ is stable, then $V=H^0(E)$ admits an
 $\omega$-balanced basis.
\end{lemma}

\noindent
{\bf Sketch of the proof.\ }
Set $W=H^0(\det
E)$ and $\W = \Hom(\Lambda^r V, W)$.  Fix an arbitrary Hermitian
metric $h$ on $E$ and consider the $L^2$--scalar product $\langle \, ,
\,\rangle=\langle \, , \,\rangle_{h,\omega}$ on $V$ built from
$\omega$ and $h$. Let $\basis{s}=\{s_1,\ldots{},s_N\}$ be an
orthonormal basis with respect to this product and
$\phi=\phi_{\basis{s}}: M \ra \G(r,N)$ the corresponding map to
the Grassmannian. On the line bundle $\det E$ consider the metric
$k_\basis{s}$ (see \eqref{kFS})
 and let $\scalar_{\W}$ and $||\cdot ||_{\W}$ be
respectively the Hermitian inner product and the norm gotten on $\W$
using $\scalar $ on $V$ and the $L^2$--metric $\scalar_{L^2}$ on $W$.
Since $\basis{s}$ is an orthonormal basis this means that for
$\alpha \in \W$
\begin{equation}
\label{onbasis}
  ||\alpha||^2_\W = \sum_I ||\alpha(s_{i_1}, \ldots{}, s_{i_r}) ||^2 _{L^2}
\end{equation}
the sum being taken over all  $r$-indices $I=(i_1, \ldots{},
i_r)$ such that $i_1<\cdots <i_r$.  An application of the Kempf-Ness theorem \ref {kempf-ness}
ensures that the function $ ||g^{-1} \cdot T_E||_\W = \exp \bigl (
\nu(g)\bigr)$ admits a minimum on $X=\SL(V)/\SU(V)$ and is proper on
geodesics transversal to $\exp(i\, \mathfrak{k})$, where 
$\mathfrak{k}$ is the Lie algebra of the
stabilizer of $T_E$ inside $\SU(V)$. Here $\nu$ is the
Kempf-Ness function based at the point $T_E$ for the action of
$\SL(V)$ on $\projective(\W)$. By \eqref{onbasis}
\begin{equation*}
  \exp\bigl(\nu(g)\bigr) = \sum_I ||g^{-1}T_E(s_I)||_{L^2}^2 .
\end{equation*}
So for any $g\in \SL(V)$ there is some multi-index $I_g$ such that
\begin{equation*}
 ||g^{-1}T_E(s_{I_g})||_{L^2}^2  \geq \frac{  \exp\bigl(\nu(g)\bigr)}{K}
\end{equation*}
where $K=\binom{N}{r}$. If we set $\eps(g)= \sqrt{K
  \exp\bigl(-\nu(g)\bigr)}$, then
\begin{equation*}
  ||\eps(g) (g^{-1}T_E) (s_{I_g}) ||^2_{L^2} \geq 1
\end{equation*}
so
Lemma 3.6 in \cite{wang-xiaowei-balance} implies that there is a $C_1\in
\real$ such that for any $g\in \SL(V)$
\begin{equation*}
  \int_M \log || \eps (g^{-1}T_E) (s_{I_g}) ||^2_{k_\basis{s}}
  \frac{\om^n}{n!}
\geq
C.
\end{equation*}
Therefore we get the inequality
\begin{equation*}
  L(g):= \int_M \biggl ( \sum_I || \eps (g^{-1}T_E) (s_{I})
  ||^2_{k_\basis{s}} \bigr )
  \frac{\om^n}{n!} \geq \nu(g) + C_2.
\end{equation*}
It follows that the function $L$ is proper on $\SL(V)/\SU(V)$ and it
must attain its minimum at some point $g$. Then $gs_1, \ldots{},gs_N $
is the desired $\om$-balanced basis. This concludes the proof.  \fine

\subsection{The proofs}

\dimostr[ of Theorem \ref{mainteor}] Let $r$ be the rank of $E$ and
$N=h^0(E)$.  Hypotheses 1. and 2.  and Lemma \ref{wang-teo}
yield the existence of an $\omega$-balanced basis
$\basis{s}=(s_1,\ldots{},s_N)$ of 
 $V=H^0(E)$.  Denote by $\phi=\phi_{\basis{s}}$
the holomorphic map $\phi: M \ra \G(r,N)$ obtained using the sections
$\basis{s}=(s_1,\ldots{}, s_N)$. Let $F: M \ra \su(N)$ be the matrix
function $F(x)=\mu_\G\bigl ( \phi(x)\bigr )$ and let $f_{jk}$ be the
entries of $F$.  The balanced condition, as rephrased in
\eqref{eq:balanced-moment}, says that
\begin{equation*}
  \int_M f_{jk} \frac{\om^n}{n!} =0.
\end{equation*}
In fact the only use of the balanced metric is to provide us with
these test functions. Using the Rayleigh principle we get
$$
\lambda_1 (M, g) \leq \frac {\int_M |\nabla f_{jk}|^2
  \frac{\omega^n}{n!}}  {\int_M |f_{jk}|^2 \frac{\omega^n}{n!}}.
$$
Thus
\begin{equation}
  \label{rayleigh-sum} \lambda_1 (M, g) \cdot \Biggl(\  \sum_{j,k=1}^N
\int_M
 |f_{jk}|^2 \frac{\omega^n}{n!} \Biggr )  \leq
  \sum_{j,k=1}^N
\int_M
 |\nabla f_{jk}|^2 \frac{\omega^n}{n!}.
\end{equation}
We claim that
\begin{equation}
  \label{phi-1}
\sum_{j,k=1}^N  \int_M     | f_{jk}|^2 \frac{\omega^n}{n!} =
\frac{r(N-r)}{N}\vol(M,g)
\end{equation}
and
\begin{equation}
  \label{nabla-phi}
  \sum_{j,k=1}^N \int_M  |\nabla f_{jk}|^2 \frac{\omega^n}{n!} =
\frac{4\pi}{ (n-1)!}\langle
  c_1(E) \cup [\omega]^{n-1}, [M]\rangle.
\end{equation}
To prove \eqref{phi-1} observe that
$$
\sum_{j,k=1}^N | f_{jk}(x)|^2 = || F(x)||^2 = || \mu_\G\bigl (\phi(
x)\bigr)||^2.
$$
Since the moment map is $ \SU(N)$-equivariant its norm is constant on
the Grassmannian.  Calculating at the point $x_0$ with Stiefel matrix
$(I_r \ 0)^t$ we get
$$
\sum_{j,k=1}^N | f_{jk}(x)|^2 = ||F(x_0)||^2=\frac{r(N-r)}{N}.
$$
From this
\eqref{phi-1}
follows immediately.

In order to prove (\ref{nabla-phi}) observe that for any $f\in
C^{\infty}(M, \complex)$ we have
\begin{equation}
\label{barbarbar}
|\nabla f|^2\om^n =
n \bigl (   i \partial f \wedge \bar\partial \bar f +i \partial
\bar{f}   \wedge \bar\partial f \bigr ) \wedge \omega^{n-1}.
\end{equation}
Therefore
\begin{gather*}
 |\nabla f_{jk}|^2 \frac{\omega^n}{n!}
 =\frac{1}{(n-1)!}
\bigl ( i \partial f_{jk} \wedge \bar\partial
\bar  f_{jk}
+ i
 \partial \bar f_{jk} \wedge \bar\partial
  f_{jk} \bigr) \wedge  \omega^{n-1}=\\
=-\frac{i}{(n-1)!} \phi^* \bigl ( \partial \mu_{jk}\wedge\debar{}\mu_{kj} +
\partial \mu_{kj}\wedge\debar \mu_{jk} \bigr ) \wedge \om^{n-1}.
\end{gather*}
Using \eqref{omG-mu} we get
\begin{gather*}
  \sum_{j,k=1}^N  |\nabla f_{jk}|^2 \frac{\omega^n}{n!} =
  \frac{2}{(n-1)!}\phi^*(\om_\G) \wedge \om^{n-1}.
\end{gather*}
To get \eqref{nabla-phi} it is enough to recall that $[\phi^*(\om_\G)]=2\pi
c_1(E)$.

Now substitute \eqref{phi-1} and \eqref{nabla-phi} in
\eqref{rayleigh-sum}. Recalling that $\int_M \om^n =n! \,\vol(M, g)$ one
immediately gets \eqref{mainest}.\fine

\dimostr[ of Theorem \ref{teorgrass}] \label{grass} 

\noindent
The proof will follow applying Theorem 
\ref{mainteor} to the bundle $E=U^*$ and $L=-K$ where $U=U_{r, N}$
and $K$ are respectively the universal subbundle 
and the canonical bundle on $M=\G (r, N)$.
Observe that if
$H^0(U^*)=V$ then $H^0(\det
U^*)=\Lambda^rV$. 
Therefore,  it is easily seen that 
the Gieseker point $T_{U^*}$ is simply the identity
map $Id$ of $\Lambda^r V$, the action of $a\in \SL (V)$ on $V$ is
just the pull-back and the action of $\SL (V)$ on $\Hom(\Lambda^r V,
\Lambda^r V)$ is given by
\begin{equation}
  (a\cdot \Phi) (s_1 \wedge \cdots \wedge s_r ) =\Phi (as_1\wedge \cdots{}\wedge a s_r)
\end{equation}
where $a\in\SL(V), s_j\in V$ and $\Phi\in
\Hom(\Lambda^r V, \Lambda^r V)$.

To apply Theorem \ref{mainteor} we need to check that $T_{U^*}=Id$ is
stable, i.e. its stabiliser is finite and its $\SL (V)$-orbit is
closed.

If
$a\cdot I=I$ then $as_1\wedge \cdots{}\wedge a  s_r=s_1\wedge \cdots{}
\wedge s_r$ for any $s_1\wedge \cdots{} \wedge s_r \in \Lambda^rV$.
It  follows that
$a=I$. Therefore the stabiliser of $T_{U^*}=I$ is trivial.

If the orbit of $\SL(V)$ through $I$ were not closed, by the 
Hilbert-Mumford criterion (see e.g.
Theorem 4.2 in \cite{birkes}) there would be a non-trivial algebraic
one-parameter subgroup $ \lambda : {\complex}^* \rightarrow \SL (V)$
such that
\begin{equation}
  \lim_{t\rightarrow 0} \lambda(t)\cdot I = T_\infty
\end{equation}
for some $T_\infty \in \Hom(\Lambda^rV, \Lambda^rV)$.  Let $s_1,
..., s_ N$ be a basis of $V$ such that $\lambda(t)s_j = t^{m_j} s_j$.
Since $\lambda$ is a 1-parameter subgroup in $\SL(V)$,
$m_1+\cdots+m_N=0$.  Assume $m_1\geq m_2 \geq ... \geq m_N$.  As
$\lambda$ is non-trivial, we have $m_1>0 > m_N$.  We claim that
$m_1+\cdots +m_r >0$. In fact, assume that $m_j\geq 0$ for $j\leq s$
and $m_{j}<0$ for $j>s$. If $s\geq r$, the sum $m_1+\cdots +m_r$ is
clearly positive, since the first term is positive and the others are
nonnegative. If instead $s<r$, then $m_{r+1} +\cdots +m_N<0$ since all
terms are negative.  Therefore
\begin{equation}
  m_1+ \cdots +m_r = - (m_{r+1} +\cdots +m_N) >0.
\end{equation}
Then we have indeed $m_1+ \cdots +m_r>0$. But then
\begin{eqnarray}
  T_\infty (s_1, \ldots{}, s_r)  = \lim_{t\rightarrow 0} \bigl (\lambda(t)
  \cdot I \bigr ) (s_1, \ldots{}, s_r) = \\
  =\lim_{t\rightarrow 0} \lambda(t)s_1 \wedge \ldots{} \wedge
  \lambda(t) s_r  =  \lim_{t\rightarrow 0} t^{m_1 +\cdots + m_k} s_1
  \wedge \ldots{} \wedge s_r.
\end{eqnarray}
This is impossible since the right hand side diverges. Therefore the
orbit is closed, and $T_{U^*}$ is stable.
The assumptions of Theorem \ref{mainteor} are
therefore satisfied.
Therefore
(\ref{mainest2}) with $\omega_g\in 2\pi c_1(M)=2\pi c_1(-K)$ and $N_r=r(N-r)=\dim M$
yields the estimate 
\begin{eqnarray*}
  \lambda_1(M,g) &\leq & 2\frac{N_r\, h^0(U^*)}{r(h^0(U^*) -r)}
\frac{\deg (U^*)} {c_1(-K)^{N_r}}
  =2\frac{N_r N}{N_r} \frac{c_1({\mathcal{O}}(1))
    \cdot c_1(-K)^{N_r-1}} {c_1(-K)^{N_r}}\\
  &=&2 \frac{c_1(-K)
    \cdot c_1(-K)^{N_r-1}} {c_1(-K)^{N_r}}=2,\\
\end{eqnarray*}
where the second equality follows from
$K= {\mathcal{O}}(-N)$.
\fdim

\section{Final remarks}

Let us indicate some lines of future research which we feel are worth
pursuing in
light of our result.

\medskip

If $M$ is a Fano manifold and $g$ is a \KE metric, then
$\lambda_1(M,g)$ $\geq 2$ and equality holds if and only if $M$ admits
nonzero holomorphic vector fields. This follows from work of Futaki,
see \cite{futaki-libro}, p.40ff.  Therefore Theorem \ref{mainteor} could
be used to rule out the existence of \KE metrics on Fano manifolds. In
fact if $M$ is an $n$-dimensional Fano manifold and $E$ a globally
generated rank $r$ vector bundle over $M$ such that
\begin{equation}\label{contrario}
  \frac{n\ h^0(E) c_1(E)\cdot c_1(-K))^{n-1}}
  {r(h^0(E) - r)c_1(-K)^{n}} < 1 ,
\end{equation}
then either the Gieseker point $T_E$ is not stable or $M$ does not
admit a \K -Einstein metric.  If $M$ does not have any nontrivial
holomorphic vector field, the equality in \eqref{contrario} is enough
to get the conclusion.

We do not know any example of a Fano manifold with a \emph{line}
bundle $E$ for which \eqref{contrario} holds. We believe such
examples, if any, would be quite interesting in view of the connection
with \KE metrics.  If $Pic(M)=\Z$ one can rule out the existence of
such line bundles using a classical result of Kobayashi and Ochiai,
according to which the index of a Fano manifold cannot exceed $n+1$.

We believe the extension to higher rank vector bundles
should on the contrary forbid some Fano manifold to have a \K
-Einstein metric.

\medskip

If $M$ is a surface of genus $g$, Yang and Yau proved that
\begin{equation}
\label{yy}
  \lambda_1(M,g)\cdot \vol(M,g)\leq 8\pi \biggl [ \frac{g+3}{2}\biggr].
\end{equation}
Optimal estimates are only known for $g=0$ or $g=1$ and for $g=1$ the
estimate \eqref{yy} is not sharp, see \cite{na}. It is therefore
natural to tackle this problem with the help of Theorem
\ref{mainteor}. Unfortunately it is not easy to construct bundles that
improve  \eqref{yy}. To 
get the best possible estimate, one has to minimize the ratio
\begin{equation*}
  \frac{h^0(E)}{r(h^0(E)-r)}\cdot \deg E 
\end{equation*}
among globally generated rank $r$ vector bundles with stable Gieseker
point. For rank $r=1$ the best possible choice is $h^0(E)=2$ and
\begin{equation*}
  \deg E =  \biggl [ \frac{g+3}{2}\biggr].
\end{equation*}
For higher rank the existence of globally generated stable bundles of
given rank and degree with fixed $h^0$ is still unanswered (see
\cite{ma1}).  But it seems very hard to improve the estimate using
vector bundles of higher rank.

\small{}

\vskip.3cm

\noindent Universit\`a di Parma,\\
\textit{E-mail:} \texttt{claudio.arezzo@unipr.it}

\vskip.3cm

\noindent Universit\`a di Milano Bicocca,\\
\textit{E-mail:} \texttt{alessandro.ghigi@unimib.it}

\vskip.3cm

\noindent Universit\`a di Cagliari,\\
\textit{E-mail:} \texttt{loi@unica.it}

\end{document}